\documentclass[12pt]{amsart}

\usepackage{amscd,amsmath,amssymb,amsfonts}

\theoremstyle{plain}

\newtheorem{thm}{Theorem}[section]
\newtheorem{theorem}[thm]{Theorem}
\newtheorem{lemma}[thm]{Lemma}

\newtheorem{proposition}[thm]{Proposition}
\theoremstyle{definition}

\newtheorem{remark}[thm]{Remark}

\numberwithin{equation}{thm}

\newcommand{\sC}{{\mathcal C}}

\newcommand{\sH}{{\mathcal H}}

\newcommand{\sK}{{\times}}

\newcommand{\sU}{{\mathcal U}}

\newcommand{\A}{{\mathbb A}}

\renewcommand{\P}{{\mathbb P}}

\newcommand{\ga}[2]{\begin{gather}\label{#1}#2 \end{gather}}

\begin{document}
\title{Short equations for the genus 2 covers of degree 3 of an elliptic curve}
\author{Jan Christian Rohde}
\address{Universit\"at Duisburg-Essen, Fachbereich Mathematik, 45117 Essen, Germany}
\email{jan.rohde@uni-essen.de}

\maketitle

\begin{abstract}
E. Kani \cite{kani} has shown that the Hurwitz functor
$\sH_{E/K,3}$, which parameterizes the (normalized) genus 2 covers
of degree 3 of one elliptic curve $E$ over a field $K$, is
representable. In this paper the moduli scheme $H_{E/k,3}$ and the
universal family are explicitly calculated over an algebraically
closed field $k$ and described by short equations.
\end{abstract}

\section*{Introduction}
T. Shaska \cite{shaska2} has given a long equation, which
describes the genus 2 curves with covers of degree 3 onto elliptic
curves over an
algebraically closed field of characteristic 0.\\
One can obtain simpler equations for the genus 2 covers of degree
3 of one elliptic curve. We use the theoretical framework, which
has been introduced by E. Kani \cite{kani}. Let $K$ be a field
with $char(K) \neq 2, 3$. We consider a $K$-scheme $S$, an
elliptic curve $E/K$ with zero point 0 and a relative genus 2
curve $C/S$. A normalized genus 2 cover $f: C \to E \sK S = E_S$
of degree 3 is a morphism of $S$-schemes of degree 3 such that
$f_*(W_{C/S}) = 2 \cdot [0_{E/S}] + E_S[2]$, where $W_{C/S}$ is
the divisor of Weierstrass points and $[0_{E/S}]$ is the zero
section of $E \times S \to S$. Two genus 2 covers $f_i : C_i \to
E_S$ are isomorphic if there is an isomorphism $\varphi : C_1 \to
C_2$ such that $f_1 = f_2 \circ \varphi$. One denotes by
$\sH_{E/K,3}(S)$ the set of isomorphism classes of normalized
genus 2 covers of degree 3 onto $E_S$. The assignment $S \to
\sH_{E/K,3}(S)$ yields a Hurwitz functor $\sH_{E/K,3}:
\underline{Sch}_{/K} \to \underline{Sets},$ which is represented
by a smooth and geometrically connected
modular curve $H_{E/K,3}$ (\cite{kani}, Theorem 1.1).\\
Let $k$ be an algebraically closed field with $char(k) \neq 2,3$.
Here $H_{E/k,3}$ and the universal family $\sC \to H_{E/k,3}$ will
be explicitly calculated. We construct and parameterize genus 2
covers of elliptic curves by suitable coverings $u : \P^1 \to
\P^1$ (Frey-Kani coverings) in Section 1. Section 2 treats the
positions of the ramification points of the Frey-Kani coverings.
The pattern
of these positions allows us to determine $H_{E/k,3}$.\\
Gerhard Frey gave some hints for the presentation of this paper.
Martin M\"oller helped to improve this paper. I would like to
thank them, and Eckart Viehweg for all his time and effort spent
in guiding me for my Diplomarbeit ("master thesis"), from which
this paper is originated.

\section{Construction of genus 2 covers}

Let us fix some elliptic curve $E$, which is given by the 4
different points $0, 1, \lambda, \infty \in \mathbb{P}^1$, over an
algebraically closed field $k$ with $char(k) \neq 2,3$. We
consider a normalized genus 2 cover $f : C \to E$ of degree 3.
There exist covers $h:C \to \P^1$ and $i:E\to \P^1$ of degree 2
and there is a cover $v: \P^1 \to \P^1$ of degree 3, which is
called "Frey-Kani covering", (see \cite{freykani}, \cite{kuhn} and
\cite{shaska}) such that this diagram commutes: \ga{2.1*}{
\begin{CD}
C @>f>> E\\
 @V h VV @V i VV\\
\P^1 @>v>> \P^1
\end{CD}
}
 Let the zero point of $E$ lie over $\infty$, and 0, 1
and $\infty$ be ramification points of a cover $u : \P^1 \to \P^1$
of degree 3. Assume that $u(0) = 0$, $u(1) = 1$ and $u(\infty) =
\lambda$. Later we will see that almost all covers $u$ are
Frey-Kani coverings. Let us first study $u$. Then the results of the 
following calculation will be used
for a calculation and parametrization of the normalized genus 2
covers $f: C \to E$ of degree 3. We have for some $p_1, c \in k$:
\begin{equation} \label{defu}
u(x_0:x_1) = (g(x_0,x_1):c(x_1-p_1x_0)x_1^2)\end{equation} with
\begin{equation*} g(x_0,x_1) = x_1^3 +
g_2x_1^2x_0+g_1x_1x_0^2+g_0x_0^3
\end{equation*}
Therefore we get by $u(0:1) = (1:\lambda)$, which implies that
$\lambda = c$, and by $u(1:1) = (1:1)$:
\begin{equation} \label{p4} p_1 =
1-\frac{g(1,1)}{\lambda}
\end{equation}

Now we want to consider the situation over the branch points 1 and
$\lambda$ of $u$. Therefore we define $\grave{u} := (x_0:x_1-x_0)
\circ u$ and $\hat{u} := (x_0:x_1-{\lambda}x_0) \circ u$. Then one
has $\grave{u}(0) = -1$, $\grave{u}(1) = 0$, and $\grave{u}(\infty) =
\lambda-1$. We obtain for some $p_2, d \in k$:
\begin{equation*}
\grave{u}(x_0:x_1) = (g(x_0,x_1):d(x_1-p_2x_0)(x_1-x_0)^2)
\end{equation*}
Thus, we get by $\grave{u}(0:1) = (1:\lambda-1)$, which implies
$\lambda-1 = d$, and by $\grave{u}(1:0) = (1:-1)$:
\begin{equation} \label{p5}
p_2 = \frac{g_0}{\lambda-1}
\end{equation}
One has $\hat{u}(0) = -\lambda$, $\hat{u}(1) = 1-\lambda$, and
$\hat{u}(\infty) = 0$. We obtain for some $e, p_3 \in k$:
\begin{equation*}
\hat{u}(x_0:x_1) = (g(x_0,x_1):e(p_3x_1-x_0)x_0^2)
\end{equation*}
Therefore we get by $\hat{u}(1:0) = (1:-\lambda)$, which implies
that ${\lambda}g_0 = e$, and by $\hat{u}(1:1) = (1:1-\lambda)$:
\begin{equation} \label{p6}
 p_3 =
1+\frac{(1-\lambda)g(1,1)}{{\lambda}g_0}
\end{equation}

By the definitions of $\grave{u}$ and $\hat{u}$, and the preceding
results, we get:
\begin{equation*}
(\lambda-1)(x_1-\frac{g_0}{\lambda-1}x_0)(x_1-x_0)^2 =
\lambda(x_1-(1-\frac{g(1,1)}{\lambda})x_0)x_1^2-g(x_0,x_1),
\end{equation*}
\begin{equation*}
{\lambda}g_0((1+\frac{(1-\lambda)g(1,1)}{{\lambda}g_0})x_1-x_0)x_0^2
=
\lambda(x_1-(1-\frac{g(1,1)}{\lambda})x_0)x_1^2-{\lambda}g(x_0,x_1)
\end{equation*} These equations of polynomials imply the following equations of coefficients of $x_1^2x_0$:
\begin{equation} \label{ustrich1}
-2(\lambda-1)- g_0 = g(1,1)-\lambda-g_2,
\end{equation}
\begin{equation} \label{baru1}
g(1,1)-\lambda-{\lambda}g_2 = 0
\end{equation}
By $(\ref{p4})$ and $(\ref{p5})$, one has the equations $g(1,1) =
\lambda - \lambda p_1$ and $g_0 = \lambda p_2-p_2$. We substitute
for $g(1,1)$ and $g_0$ in $(\ref{ustrich1})$ and $(\ref{baru1})$
and obtain:
\begin{equation} \label{pvier}
p_1 = -g_2 = -2(\lambda-1)+p_2-{\lambda}p_2+p_1\lambda
\end{equation}
\begin{equation*}
\Rightarrow(\lambda-1)(-2-p_2+p_1) = 0
\end{equation*}
One can divide by $\lambda-1$, because we have $\lambda \neq 1$.
Thus, we have:
 \begin{equation}\label{suessegerade} -2-p_2+p_1 = 0 \Leftrightarrow
p_2 = p_1-2
\end{equation}
One obtains by ($\ref{p4}$), ($\ref{p5}$) and ($\ref{p6}$):
\begin{equation} \label{nochwas} p_3 = 1+\frac{(1-\lambda)g(1,1)}{{\lambda}g_0} =
1-\frac{1-p_1}{p_1-2} = \frac{2p_1-3}{p_1-2} \end{equation}

The equations $p_2 = p_1-2$ and $p_2 = \frac{g_0}{1-\lambda}
\Leftrightarrow p_2(\lambda-1) = g_0$ imply:
\begin{equation} \label{g0} ({\lambda}-1)(p_1-2) = g_0
\end{equation}
Using the equations ($\ref{p4}$), ($\ref{pvier}$) and ($\ref{g0}$)
we get:
\begin{equation} \label{g1}
g_1 = g(1,1) - 1 -g_2-g_0 = \lambda-p_1\lambda - 1 + p_1 +
(1-{\lambda})(p_1-2) =-2p_1\lambda + 2p_1 + 3\lambda-3
\end{equation}
Thus, by ($\ref{pvier}$), ($\ref{g0}$), and ($\ref{g1}$), we obtain
\begin{equation} \label{geh}
g(1,x) = x^3-p_1x^2+(-2{\lambda}p_1+2p_1+3\lambda-3)x-(1-{\lambda})(p_1-2).
\end{equation}
Hence the cover $u$ is completely determined by ($\ref{defu}$), $\lambda
= c$, and ($\ref{geh}$). Now we can
apply the results of this calculation and construct normalized
genus 2 covers of degree 3. Later we will use that we get by ($\ref{geh}$):
\begin{equation} \label{babu}
\lambda(x-p_1)x^2 - \lambda g(1,x) =
-\lambda(1-\lambda)((2p_1-3)x-p_1+2),\end{equation}
\begin{equation}\label{bubu} \lambda(x-p_1)x^2 - g(1,x) = (\lambda-1)(x-p_1+2)(x-1)^2
\end{equation}
Now we consider smooth curves $C$ of genus 2, which are given by
\begin{equation} \label{faber} y^2 = (x-p_1)(x-p_1+2)((2p_1-3)x-p_1+2)g(1,x) \end{equation}
for some $p_1$.
\begin{proposition} \label{ututs}
The curve $C$ given in ($\ref{faber}$) is smooth if $u$ is not
ramified over $\infty$. Let $h:C \to \P^1$ and $i:E \to \P^1$ be
the natural projections, which are given by $(x,y) \to x$. Then
the normalized covers $f_{\pm}:C\to E$ of degree 3 with Frey-Kani
covering $u$ for $h$ and $i$ are given by
\begin{equation*} f_{\pm}(x,y) = (g(1,x):
\lambda(x-p_1)x^2:\pm\frac{\lambda(\lambda -1)yx(x-1)}{g(1,x)}).
\end{equation*}
\end{proposition}
\begin{proof} Using that $u$ is unramified over $\infty$ we see that
$g(1,x)$ has 3 different zeros.
By $(\ref{suessegerade})$ and $(\ref{nochwas})$, we conclude that
$u$ maps the 3 other points in $\P^1$, which give the Weierstrass points of
$C$, to 3 different points. Therefore the Weierstrass points of $C$ are given
by 6 different points of $\P^1$, and $C$ is smooth. Let $X_1 =
\lambda(x-p_1)x^2$ and $X_0 = g(1,x)$. Note that $E = 
\{y^2x_0 - x_1(x_1-x_0)(x_1-\lambda x_0) = 0\} \subset \P^2$. There is a cover $f: C \to
E$ with Frey-Kani covering $u$ for the chosen degree 2 covers $h$
and $i$ if and only if there is an $Y\in k(C)$, which satisfies
the equation $Y^2X_0 = X_1(X_1-X_0)(X_1-\lambda X_0)$ in $k(C)$.
By ($\ref{babu}$) and ($\ref{bubu}$), one can easily check that

this equation has the solutions:
\begin{equation*} Y = \pm \frac{\lambda(\lambda -1)yx(x-1)}{g(1,x)}
\end{equation*}
By the definitions of $C$ and $f_{\pm}$, one can easily see that
$f_+$ and $f_{-}$ map the 3 different Weierstrass points of $C$,
which are given by the zeros of $g(1,x)$, to $(0:0:1) = 0_E$.
Thus, $f_+$ and $f_{-}$ are normalized.
\end{proof}
Let $\sU$ be the set of coverings $\P^1 \to \P^1$ of degree 3,
which satisfy our choice of coordinates of ramification points and
do not have any ramification point over $\infty$.

\begin{proposition} \label{klasse}
Using Proposition $\ref{ututs}$ we have a bijection $\chi : \sU
\to \sH_{E/k,3}(k)$. This map is given by $u \to [f_{\pm}:C \to
E]$.
\end{proposition}
\begin{proof}
By the hyperelliptic involution on $C$, we conclude that the
covers $f_{\pm}:C \to E$ lie in the same isomorphism class, and
the map $\chi : \sU \to \sH_{E/k,3}(k)$ is well defined.\\ Let
$f:C \to E$ be a normalized genus 2 cover of degree 3. The
Frey-Kani covering of $f$ is not ramified over $\infty$ and has
ramification points over 0, 1 and $\lambda$ (see \cite{frey}). One
can put these ramification points on 0, 1 and $\infty$ such that
the Frey-Kani covering is some $u \in \sU$.
In \cite{frey}, page 92-93 the curve $C$ is described by the
Frey-Kani covering. This description of
$C$, ($\ref{babu}$) and ($\ref{bubu}$) imply that $C$ is given by
($\ref{faber}$) for some $p_1$. One can assume that $h$ and $i$
are given by $(x,y) \to x$. Hence by Proposition $\ref{ututs}$, we
conclude that $\chi$ is
surjective.\\
Let $\chi(u_1) = [f_1:C_1 \to E] = [f_2:C_2 \to E] = \chi(u_2)$
and $h_i : C_i \to \P^1$ be the natural degree 2 covers (for $i =
1,2$). Then there is an isomorphism $i:C_1 \to C_2$ such that $f_1
= f_2 \circ i$ and an $a \in Aut(\P^1)$ such that $a \circ h_1 =
h_2 \circ i$ (see \cite{hart}, page 304, IV, Exercise 2.2.(a)).
Therefore we conclude $u_1 = u_2 \circ a$. Let $x \in\{0, 1,
\infty\}$. We get by our assumptions that $u_1(x) = u_2(x)$. These
points are ramification points of $u_1$ and $u_2$. Thus, $u_1 =
u_2 \circ a$ implies that $a(x) = x$ ($\forall\mbox{ } x \in\{0,
1, \infty\}$). Therefore we have $a = \rm{id}$ and $\chi$ is
injective.
\end{proof}

\section{The reckoning of the moduli space}

 By $u(\infty) = \lambda \neq 0 =
u(p_1)$, we conclude that $p_1$ and $\infty$ can not coincide. Thus,
($\ref{defu}$), $c= \lambda$, and ($\ref{geh}$) imply that $g(1,x)$ is determined by $p_1$ and
$\lambda$, and that:

\begin{lemma} \label{univf}
We have an injective map $\imath: \sU \to\A^1$, which maps the
morphism, which is given by $x \to
\frac{\lambda(x-p_1)x^2}{g(1,x)}$, to $p_1$.
\end{lemma}

\begin{remark} \label{hahar}
The morphism, which is given by $x \to
\frac{\lambda(x-p_1)x^2}{g(1,x)}$ for some $p_1$, has the degree 3
if and only if we have $g(1,p_1) \neq 0$ and $g(1,0) \neq 0$. By 
($\ref{geh}$), one can easily see
that this is true if and only if $p_1 \in k \setminus \{1,2\}$.
Thus, by ($\ref{babu}$) and ($\ref{bubu}$), we have $p_1 \in
\imath(\sU)$ if and only if we have $p_1 \in k \setminus \{1,2\}$
and $g(1,x)$ has 3 different zeros.
\end{remark}

Let $p_1 \neq 1, 2, \infty$, and $u$ be the morphism, which is
given by $x \to \frac{\lambda(x-p_1)x^2}{g(1,x)}$. The preceding
remark, ($\ref{babu}$) and ($\ref{bubu}$) imply that $u$ has the
degree 3 and satisfies our choice of coordinates of ramification
points, and that we have $p_1 \in \imath(\sU)$ if and only if $u$
does not have a ramification point over $\infty$. By the Hurwitz
formula, $u$ has a ramification divisor of degree 4. But we do not
know the position of a fourth branch point. Let $u(\mu) =
u(\delta) = \zeta$ and $\delta$ be a ramification point of $u$. We
have $(x_0:x_1-{\zeta}x_0) \circ u = (g:F(x_1 - \mu x_0) (x_1 -
{\delta}x_0)^2)$ for some $F \in k$ resp.,
\begin{equation*}
F(x_1-{\mu}x_0)(x_1-{\delta}x_0)^2 =
(\lambda-\zeta)x_1^3+(g(1,1)-\lambda-\zeta g_2) x_1^2x_0-\zeta
g_1x_1x_0^2-{\zeta}g_0x_0^3.
\end{equation*}

This equation of polynomials implies the following equations of
coefficients:
\begin{eqnarray*}
F = \lambda-\zeta\\
 -F\mu-2F\delta=
g(1,1)-\lambda-\zeta g_2\\
 F\delta^2+2F\delta\mu =
-\zeta g_1\\
 -\mu\delta^2F= -{\zeta}g_0
\end{eqnarray*}
Using $F = \lambda-\zeta$, $p_1 = -g_2$, $-p_1\lambda =
g(1,1) - \lambda$, ($\ref{pvier}$), and ($\ref{p4}$) we substitute in the
second equation:
\begin{equation*}
(-\lambda+\zeta)\mu+2(-\lambda+\zeta)\delta=
-p_1\lambda+\zeta{p_1}
\end{equation*}
Let $\lambda\neq\zeta$. Then we have:
\begin{equation} \label{mue}
\mu = p_1-2\delta
\end{equation}
By $\mu = p_4-2\delta$ and and the same substitution as above, we obtain:
\begin{equation} \label{1gl}
(\lambda-\zeta)(2\delta{p_1}-3\delta^2) =
-\zeta(1-\lambda)(2p_1-3),
\end{equation}
\begin{equation} \label{2gl}
(2\delta^3-p_1\delta^2)(\lambda-\zeta)= {\zeta}(1-\lambda)(p_1-2)
\end{equation}
By $(\ref{2gl})$, we get the equation $\zeta =
\frac{(2\delta^3-p_1\delta^2)(\lambda-\zeta)}{(1-\lambda)(p_1-2)}$.
Thus, one has by ($\ref{nochwas}$) and ($\ref{1gl}$):
\begin{equation*}
(\lambda-\zeta)(2\delta{p_1}-3\delta^2) =
(2\delta^3-p_1\delta^2)(\lambda-\zeta)(-p_3)
\end{equation*}
The solution $\zeta = \lambda$ gives the ramification point
$\infty$ and $\delta = 0$ gives the ramification point 0.
Therefore we can divide by $\lambda-\zeta$ and $\delta$, and get
\begin{equation*}
2{p_1}-3\delta = (2\delta^2-p_1\delta)(-p_3) \Leftrightarrow 0 =
\delta^2-\frac{p_1p_3+3}{2p_3}\delta+\frac{p_1}{p_3} =
(\delta-1)(\delta-\frac{p_1}{p_3}),
\end{equation*}
because we get by ($\ref{nochwas}$):
\begin{equation*}
\frac{p_1p_3+3}{2p_3} = \frac{p_1(2p_1-3)+3(p_1-2)}{2(2p_1-3)} =
\frac{(2p_1^2-4p_1)+(4p_1-6)}{2(2p_1-3)} = \frac{p_1}{p_3}+1
\end{equation*}
Thus, a ramification point $\delta$ of $u$ is given by
\begin{equation*}
\delta = \frac{p_1}{p_3} = \frac{p_1(p_1-2)}{2p_1-3}.
\end{equation*}

We substitute $\delta = \frac{p_1(p_1-2)}{2p_1-3}$ in
($\ref{1gl}$) and get:
\begin{equation*}
(\lambda-\zeta)(2\frac{p_1(p_1-2)}{2p_1-3}p_1-3(\frac{p_1(p_1-2)}{2p_1-3})^2)
= -\zeta(1-\lambda)(2p_1-3)\end{equation*}
\begin{equation*}\Leftrightarrow (\lambda-\zeta)p_1^3(p_1-2) =
-\zeta(1-\lambda)(2p_1-3)^3
 \Leftrightarrow {\lambda}p_1^3(p_1-2)
= \zeta((\lambda-1)(2p_1-3)^3+p_1^3(p_1-2))\end{equation*}
\begin{equation}\label{degen}\Leftrightarrow
\zeta =
\frac{{\lambda}p_1^3(p_1-2)}{(\lambda-1)(2p_1-3)^3+p_1^3(p_1-2)}
\end{equation}

Recall that $E$ is
determined by $\lambda$, and that $E$ resp., $\lambda$ is fixed.
Hence by $\ref{degen}$, the set of values of $p_1$, which induce a 4th ramification
point over $\infty$, is given by $Z = \{(\lambda-1)(2\cdot p_1 -3)^3+p_1^3(p_1-2)= 0\} \subset \P^1$. 

\begin{remark} \label{4losungen}
Four different points lie in $Z$.
\end{remark}
\begin{proof}
We consider the morphism $\tilde{\zeta}:\P^1 \to \P^1$, which is
given by
\begin{equation*}
\tilde{\zeta}(p_1) =
\frac{{\lambda}p_1^3(p_1-2)}{(\lambda-1)(2p_1-3)^3+p_1^3(p_1-2)} =
(1 +
\frac{\lambda-1}{\lambda}\frac{(2p_1-3)^3}{p_1^3(p_1-2)})^{-1}.
\end{equation*}
Therefore the ramification points of $\tilde{\zeta}$ are the
ramification points of the morphism $\sigma$, which is given by
$p_1 \to \frac{(2p_1-3)^3}{p_1^3(p_1-2)}$. The derivative of
$\sigma$ is
\begin{equation*}
\sigma^{\prime}(p_1) = \frac{6(2p_1-3)^2p_1^3(p_1-2)
-(2p_1-3)^3(3p_1^2(p_1-2) + p_1^3)}{p_1^6(p_1-2)^2} =
\frac{-2(2p_1-3)^2p_1^2(p_1 - 3)^2}{p_1^6(p_1-2)^2}.
\end{equation*} Thus, the ramification points of $\tilde{\zeta}$ are $0$,
$\frac{3}{2}$ and $3$. The statement follows by the fact that all
$x \in \{0, \frac{3}{2}, 3\}$ fulfil $(\lambda-1)(2\cdot x
-3)^3+x^3(x-2) \neq 0$.
\end{proof}

 By Lemma $\ref{univf}$ and
Remark $\ref{hahar}$, we have a bijection between $\sU$ and the
set of closed points of $\P^1 \setminus (\{1,2, \infty\} \cup Z)$,
where $\P^1 \setminus (\{1,2, \infty\} \cup Z)$ is an (Zariski)
open subset of $\P^1$. Recall that $\sH_{E/k,3}$ is represented by
a smooth modular curve $H_{E/k,3}$ (see \cite{kani}, Theorem 1.1).
Therefore, we conclude by Proposition $\ref{ututs}$ and
Proposition $\ref{klasse}$:

\begin{theorem} \label{tutu}
We have ${H}_{E/k,3} \cong \P^1 \setminus (\{1,2, \infty\} \cup
Z) $. The fiber $\sC_{p_1}$ of the universal family $\sC \to \P^1 \setminus (\{1,2, \infty\} \cup
Z)$ is given by ($\ref{faber}$) for all $p_1 \in \P^1 \setminus (\{1,2, \infty\} \cup
Z)$.
\end{theorem}

\end{document}